\documentclass[a4paper,10pt]{article}
\usepackage{ifstyle}
\usepackage{graphicx}
\usepackage{epsfig} 
\numberwithin{equation}{section} 

\usepackage[english]{babel}               
\usepackage[T1]{fontenc}
\usepackage[latin1]{inputenc}
\usepackage{amssymb}   
\usepackage{array}                 
\def\XXint#1#2#3{{\setbox0=\hbox{$#1{#2#3}{\int}$}
\vcenter{\hbox{$#2#3$}}\kern-.5\wd0}}

\newtheorem{theorem}{Theorem} 
\newtheorem{corol}{Corollary}

\newtheorem{lem}{Lemma}
\newtheorem{Claim}{Claim}

\def\opn#1#2{\def#1{\operatorname{#2} } } 
\opn\Rm{Rm}
\opn\Ch{Ch}
\opn\Td{Td}
\opn\Ric{Ric} 
\opn\Rc{Rc}
\opn\rk{rk}
\opn\Scal{Sc}
\opn\SB{SB}
\opn\Tr{Tr}
\opn\Trac{Tr} 
\opn\det{det} 
\opn\diam{diam} 
\opn\dist{dist} 
\opn\div{div}
\opn\dim{dim}
\opn\Ker{Ker} 
\opn\NS{NS}
\opn\loc{loc}
\opn\For{For}
\opn\Hil{Hil}
\opn\exp{exp}
\opn\mult{mult}
\opn\Psh{Psh}
\opn\Vol{Vol}
\opn\exph{exph}
\opn\Herm{Herm}
\opn\End{End} 
\opn\Hess{Hess} 
\opn\Vol{Vol} 

\newcommand{\N}{\mathbb N}
\newcommand{\Z}{\mathbb Z}  
\newcommand{\Q}{\mathbb Q}
\newcommand{\R}{\mathbb R}
\newcommand{\C}{\mathbb C}

\newcommand{\contract}{\mathrel{\kern-1.5pt\vrule width6.0pt height0.4pt depth0pt
                \vrule width0.4pt height4.0pt depth0pt}}
\newcommand{\retract}{\mathrel{\kern-1.5pt\vrule width0.4pt height4.0pt depth0pt
                          \vrule width6.0pt height0.4pt depth0pt}}

\newdimen\boxrulethickness \boxrulethickness=.07em
\newdimen\Openboxwidth \Openboxwidth=.72em
\newcommand{\Openbox}{%
 \leavevmode
 \hbox{%
   \hfil\vrule width\boxrulethickness
   \vbox to\Openboxwidth{%
     \advance\Openboxwidth -2\boxrulethickness
     \hrule height \boxrulethickness width\Openboxwidth\vfil
     \hrule height\boxrulethickness}%
   \vrule width\boxrulethickness\hfil
 }}

\begin{document}
\begin{center} 
%\Large{\bf{K\"ahler-Einstein metrics over singular K\"ahler spaces}}
%\Large{\bf{Degenerate K\"ahler-Einstein metrics over varieties of general type}}
\Large{\bf{Lecture notes on the Ein-Popa extension result}}
\\
\vspace{0.4cm}
\large{Nefton Pali}
\end{center} 
\begin{abstract}
These are lecture notes on a recent remarkable preprint of Ein-Popa, which simplifies the algebraic proof of the finite generation of the canonical ring given by the team BCHM.
The Ein-Popa extension result has been translated in the analytic language by Berndson-Paun and Paun. 
In these notes we follow the analytic language used in Berndson-Paun and Paun.
The author of this manuscript does not claim any originality of the main ideas and arguments which are due to Ein-Popa, based in their turn in the ideas of Hacon-McKernan, Takayama and Siu.
\end{abstract}
\section{The Ein-Popa extension result}
Let $X$ be a complex manifold. 
The multiplier ideal sheaf ${\cal I}(\theta)\subset {\cal O}_X$ associated to a closed positive $(1,1)$-current $\theta$ is the sheaf of germs of holomorphic functions $f\in {\cal O}_x$ such that 
$$
\int_{U_x}|f|^2e^{-\varphi}<+\infty\,,
$$
where $\varphi$ is a local potential of the current $\theta=\frac{i}{2\pi}\partial\bar\partial \varphi$ over some neighborhood $U_x$ of the point $x$. Let $Z\subset X$ be a smooth hypersurface. 
If the restriction of the local potentials of $\theta$ to $Z$ is  not identically $-\infty$ on any local connected component of $Z$
we can also define the multiplier ideal sheaf ${\cal I}(\theta_{{|Z}})\subset {\cal O}_Z$ in a similar way. In this setting we introduce the adjoint ideal sheaf 
${\cal I}_Z(\theta)\subset {\cal I}(\theta)$ of germs of holomorphic functions $f\in {\cal I}(\theta)_x$ such that 
$$
\int_{Z\cap U_x}|f|^2e^{-\varphi}<+\infty\,.
$$
We will use also the analogue notations ${\cal I}(\psi_{{|Z}})\subset {\cal O}_Z\,,\,{\cal I}_Z(\psi)\subset {\cal I}(\psi)$ with respect to a global quasi-plurisubharmonic function $\psi$ which is  not identically $-\infty$ on any connected component of $Z$. %The current $\theta$ will be called KLT if $\int_{U_x}e^{-\varphi}<+\infty$ for all local potentials $\varphi$.
We observe now the following claim.
\begin{Claim}\label{OT-Bas-ext}
Let $Z\subset X$ be a smooth and irreducible hypersurface inside a complex projective manifold $X$ and let $L$ be a line bundles over $X$ such that the class $c_1(L)$ admits a
K\"{a}hler current $\theta$ with well defined restriction $\theta_{|Z}$.
Then the restriction map 
$$
H^0\left(X\,,\,{\cal O}_X(K_X+Z+L)\otimes {\cal I}_Z(\theta)\right)\longrightarrow H^0\left(Z\,,\,{\cal O}_Z(K_Z+ L_{{|Z}})\otimes{\cal I}(\theta_{{|Z}})\right)\,,
$$
is surjective. 
\end{Claim}
This claim is a direct consequence of the adjunction formula and the singular version of the Ohsawa-Takegoshi-Manivel extension theorem \cite{Pal}, \cite{Mc-Va}, \cite{Man}, \cite{Dem2},\cite{Oh-Ta}, \cite{Oh}. 

We will note by $\lambda(\theta, A):=\inf_{x\in A}\lambda(\theta)_x$ the generic Lelong number of $\theta$ along an irreducible complex analytic set $A\subset X$. Similar notations $\lambda(\psi)_x$ and $\lambda(\psi, A)$ will be employed also for global quasi-plurisubharmonic functions.
We remind that a quasi-plurisubharmonic functions $\psi$ is called with analytic singularities 
if it can locally be expressed as 
\begin{eqnarray*}%\label{exp-an-sg}
\psi=c\log\sum_j|h_j|^2+\rho\,,
\end{eqnarray*}
with $c\in \R_{>0}$, $h_j$ holomorphic functions and $\rho$ a bounded function. A closed positive $(1,1)$-current is called with analytic singularities is its local potentials has this property.
We remind the following well known (from algebraic geometers) claim (see also \cite{Be-Pa}, \cite{Pa2} for a similar statement).
\begin{Claim}\label{cv-int-Lel}
Let  $B\subset \C^n$ be an open ball, let $V\subset B$ be a hyperplane, let $v=0$ be its equation and let $\psi$ be a  quasi-plurisubharmonic function
with analytic singularities over $B$ which is not identically $-\infty$ over $V$ and let 
$$
\Lambda_{\Omega}:=\sup_{x\in \Omega}\lambda(\psi)_x<+\infty\,,
$$
for any relatively compact open set $\Omega\subset\subset B$.
Then 
$$
I_{\varepsilon, \delta}:=\int_{\Omega} |v|^{-2(1-\varepsilon)}e^{-\delta\psi}<+\infty\,,
$$
for all $\varepsilon\in (0,1)$ and $\delta\in (0,2/\Lambda_{\Omega})$.
\end{Claim}
$Proof$. The assumptions on the restriction to $V$ of the function 
$$
\psi=c\log\sum_j|h_{j}|^2+\rho\,,
$$
implies the existence of a blow-up map $\mu:(v,\zeta)\mapsto(v,z)$ such that
$$
\psi\circ \mu=c\log|\zeta^{\alpha}|^2+R\,, 
$$
with $\alpha\in \Z^{n-1}_{\ge 0}$ and $R$ a bounded function.  This last equality follows from the fact that we can construct the blow-up map $\mu$ in a way that the sheaf
$
\mu^*\sum_j{\cal O}\cdot h_{ j}\,
$
is invertible. %($\beta_j=0\Rightarrow \alpha_j=0$ since $V({\cal J}({\varphi}))\subset \varphi^{-1}(-\infty)$) 
Moreover we can also assume that the Jacobian $J(\mu)$ of $\mu$ equal to
a monomial $\zeta^\beta$, $\beta\in \Z^{n-1}_{\ge 0}$ up to an invertible factor. On the other hand Skoda's lemma implies
$$
+\infty\;>\;\int_{\Omega}e^{-\delta\psi}\;=\;\int_{\mu^{-1}(\Omega)}|\zeta^{\alpha}|^{-2c\delta}\,|J(\mu)|^2\,e^{-\delta R}\,,
$$
for all $\delta\in (0,2/\Lambda_{\Omega})$.
Thus Fubini's formula implies that the integrability of the function 
$
f_{\delta}:=|\zeta^{\alpha}|^{-2c\delta}|\zeta^{\beta}|^2
$
is a sufficient condition for the convergence of the integrals $I_{\varepsilon, \delta}$.\hfill$\Box$
\\
\\
We prove now the following analytic version of the Ein-Popa \cite{Ei-Po} extension result (see \cite{Be-Pa} and \cite{Pa2} for similar statements). The idea of the proof is due to Ein-Popa \cite{Ei-Po}, based in their turn in the ideas of Hacon-McKernan \cite{Ha-Mc1}, \cite{Ha-Mc2}, Takayama \cite{Tak} and Siu \cite{Siu1}, \cite{Siu2}.
We will follow closely the translation in the analytic language by Berndson-Paun \cite{Be-Pa} and Paun \cite{Pa1}, \cite{Pa2}. 
%%%%%%%%%%%%%%%%%%%%%%%%%%%%%%%%%%%%%%%%%%%%%%%%%%%%%%%%%%%%%%%%%%%%%%%%%%%%%%%%%%%%%%%%%%%%%%%%%%%%%%%%%%%%%%%%%%%%%%%%%%%%%%%%%%%%%%%%%%%%%%%%%%%%%%%%%%%%%%%%%%%%%%%%%%%%%%%%%%%%%
%%%%%%%%%%%%%%%%%%%%%%%%%%%%%%%%%%%%%%%%%%%%%%%%%%%%%%%%%%%%%%%%%%%%%%%%%%%%%%%%%%%%%%%%%%%%%%%%%%%%%%%%%%%%%%%%%%%%%%%%%%%%%%%%%%%%%%%%%%%%%%%%%%%%%%%%%%%%%%%%%%%%%%%%%%%%%%%%%%%%%
%%%%%%%%%%%%%%%%%%%%%%%%%%%%%%%%%%%%%%%%%%%%%%%%%%%%%%%%%%%%%%%%%%%%%%%%%%%%%%%%%%%%%%%%%%%%%%%%%%%%%%%%%%%%%%%%%%%%%%%%%%%%%%%%%%%%%%%%%%%%%%%%%%%%%%%%%%%%%%%%%%%%%%%%%%%%%%%%%%%%%
%%%%%%%%%%%%%%%%%%%%%%%%%%%%%%%%%%%%%%%%%%%%%%%%%%%%%%%%%%%%%%%%%%%%%%%%%%%%%%%%%%%%%%%%%%%%%%%%%%%%%%%%%%%%%%%%%%%%%%%%%%%%%%%%%%%%%%%%%%%%%%%%%%%%%%%%%%%%%%%%%%%%%%%%%%%%%%%%%%%%%
\newpage
\begin{lem}\label{Hk-Mk-lem}
Let $X$ be a complex projective manifold, let $Z\subset X$ be a smooth irreducible hypersurface and let $L$ be a holomorphic $\Q$-line bundle over $X$ such that;
\\
\\
{\bf I)}
there exists a closed positive $(1,1)$-current $\Theta \in c_1(K_X+Z+L)$  with well defined restriction $\Theta_{|Z}$,
\\
\\
{\bf II)} there exists a decomposition of $\Q$-line bundles $L={\cal O}_X(\Delta)+{\cal O}_X(D)+R$, where;
\\
\\
$\bullet$ 
$
\Delta=\sum_{j=1}^N \lambda_j Z_j
$
is a divisor over $X$ with $\lambda_j\in \Q\cap [0,1)$ and $Z_j\subset X$ distinct irreducible smooth hypersurfaces with normal crossing intersection with $Z$ such that $Z\cap Z_j\cap Z_l =\emptyset$ for all $j\not = l$. 
\\
\\
$\bullet$ $D$ is an effective $\Q$-divisor over $X$ such that $Z$ is not one of its components and $R$ is a holomorphic $\Q$-line bundle over $X$ which admits a K\"{a}hler current
$
\rho\in c_1(R)
$
with bounded local potentials over $Z$.%well defined restricton $\theta_{|Z}$.% such that ${\cal I}(m^{-1}\theta_{|Z})={\cal O}_Z$,
\\
\\
Let also $(V_t)_{t=1}^{N'}$ be the irreducible components of the family $(Z_j\cap Z)_{j=1}^N$ and let $m\in \N_{>1}$ such that $m\lambda_j\in \N$ for all $j$, $mD$ is integral and $mR$ is a holomorphic line bundle.
Then  for any section 
$$
u\in H^0\Big(Z\,,\,m\left(K_Z+ L_{{|Z}}\right)\Big)\,,
$$
with the vanishing property
\begin{eqnarray}\label{Sec-div}
\div u\;-\; m\left(\,\sum_{t=1}^{N'} \lambda(\Theta_{|Z}\,,V_t)\, V_t\;+\;D_{|Z}\right)\;\ge\; 0\,,
%\div u- \sum_{t=1}^{N'} m  \mu_t V_{Z,t}\ge 0\,,\qquad \mu_t:=\max\left\{\lambda_{Z,t}\,,\,\lambda(\Theta_{|Z}, V_{Z,t})/(2\pi) \right\}\,,
\end{eqnarray}
there exists a section 
$$
U\in H^0\Big(X\,,\,m(K_X+Z+L)\Big)\,,\qquad
U_{{|Z}}=u\otimes(d\zeta)^m\,,
$$
with $\zeta\in H^0(X, {\cal O}(Z))$ such that $\div \zeta= Z$.
\end{lem}
$Proof$.
\\
{\bf Notations.} For all $\nu\in \N$ we define the integers
$$
k_{\nu}:=\max\{k\in \N\,:\,km\le \nu\}\,,
$$
and $q_{\nu}:=\nu-k_{\nu}m=0,...,m-1$. Let $\omega>0$ be a K\"{a}hler form over $X$ and set 
$$
\Omega_X:=\frac{\omega^n}{n!}\,,\qquad
\Omega_Z:=\frac{\omega^{n-1}_{{|Z}}}{(n-1)!}\,,
$$
We consider now the crucial Ein-Popa decomposition \cite{Ei-Po}
$$
{\cal O}_X(m\Delta)=L_1+\cdots + L_{m-1}\,,
$$
with $L_k:={\cal O}_X(\Delta_k)$ and with
$$
\Delta_k:=\sum_{m\lambda_j=k}Z_j\,,
$$
for all $k=1,...,m-1$. %Notice that $L_{m-1}={\cal O}_X$ by our choice of $m$.
Let $h_{Z_j}$ be smooth hermitian metrics over ${\cal O}_X(Z_j)$, let $h_{L_k}$ be the induced smooth hermitian metric over $L_k$ and let denote by $h_{L_k}e^{-\varphi_k}$ the canonical singular hermitian metric associated to the divisor $\Delta_k$. %(Notice that this is the trivial metric in the case $k=m-1$.)
%The key property of the Ein-Popa decomposition is that $\lambda(\varphi_k)_z\le 2$ for all $k=1,...,m-1$ and $z\in Z$.
We equip the line bundle
$$
L_m:={\cal O}_X(mD)+mR\,,
$$ 
with the singular hermitian metric $h_{L_m}e^{-\varphi_m}$ such that
$$
2\pi\, m\,(\,[D]+\rho)=i\,{\cal C}_{h_{L_m}}(L_m)+i\partial\bar\partial \varphi_m\,.
$$ 
(As before $h_{L_m}$ is smooth.)  Let also $h_Z$ be an arbitrary  smooth hermitian metric on ${\cal O}(Z)$. 
We equip the line bundle
$$
F_m\;:=\;m(K_X+Z)+\sum_{j=1}^m L_j\;=\;m(K_X+Z+L)\,,
$$ 
with the smooth hermitian metric 
$$
h_m:=\Omega^{-m}_X\otimes h^m_Z\otimes h_{L_1}\otimes\cdots \otimes  h_{L_m}\,.
$$
Let $(A,h_A)$ be an ample line bundle over $X$ with $0<\omega_A:=i\,{\cal C}_{h_A}(A)$ and let define  for all $\nu\in \N$ the Siu-Demailly \cite{Dem3}, \cite{Siu1}, \cite{Siu2} type line bundle 
$$
{\cal L}_{\nu}:=k_{\nu}F_m+q_{\nu}(K_X+Z)+\sum_{j=0}^{q_{\nu}}L_j\,,
$$
with $L_0:=A$. We equip the line bundle ${\cal L}_{\nu}$ with the smooth hermitian metric 
$$
H_{\nu}:= h_m^{k_{\nu}}\otimes \Omega_X^{-q_{\nu}}\otimes h^{q_{\nu}}_Z\otimes  h_{L_0}\otimes h_{L_1}\otimes \cdots \otimes h_{L_{q_{\nu}}}\,.
$$ 
We choose the line bundle $(A,h_A)$ sufficiently ample such that;
\\
\\
{\bf (A1)} for all $q=0,...,m-1$ the line bundle ${\cal L}_{q {|Z}}\equiv qK_Z+(L_0+\cdots + L_q)_{{|Z}}$ is base point free, globally generated by some family $(s_{q,j})_{j=1}^{N_q}\subset H^0(Z,{\cal L}_{q {|Z}})$,
\\
\\
{\bf (A2)} the restriction map $H^0(X,F_m+A)\longrightarrow H^0(Z,F_m+A)$ is surjective,
\\
\\
{\bf (A3)} for all $q=0,...,m-1$ hold the inequality $i\,{\cal C}_{H_q}({\cal L}_q)\ge 2\pi m\omega$.
\\
\\
We note by $|\cdot|_{\nu}$ the norm of the smooth hermitian metric 
$$
\Omega_Z^{\nu}\otimes h_{L_1}^{k_{\nu}}\otimes \cdots \otimes h_{L_m}^{k_{\nu}}\otimes h_{L_0}\otimes h_{L_1}\otimes \cdots \otimes h_{L_{q_{\nu}}}\,,
$$
over the line bundle 
$$
{\cal L}_{{\nu} {|Z}}
\;\equiv \;
\nu K_Z+k_{\nu}\sum_{j=1}^m L_{j{|Z}}+\sum_{j=0}^{q_{\nu}}L_{j {|Z}}
\;=\;
\nu K_Z+k_{\nu}mL_{{|Z}}+\sum_{j=0}^{q_{\nu}}L_{j {|Z}}\,.
$$
The assumption (A1) implies
$$
\max_{0\le p,q\le m-1}\max_Z\;\frac{\sum_{j=1}^{N_q}|s_{q,j}|^2_q}{\sum_{t=1}^{N_p}|s_{p,t}|^2_p}\;=\;C\;<\;+\infty\,.
$$
We prove now the following claim (see also \cite{Ei-Po}, \cite{Tak}, \cite{Be-Pa}, \cite{Pa1} and \cite{Pa2}).
\begin{Claim}\label{Pau-lm}
Let $u$ and $\zeta$ as in the statement of the lemma \ref{Hk-Mk-lem} and let
$$
\sigma_{\nu,j}:=u^{k_{\nu}}\otimes s_{q_{\nu},j}\in H^0(Z, {\cal L}_{{\nu} {|Z}})\,,
\qquad
j=1,...,M_{\nu}:=N_{q_{\nu}}\,.
$$
Then for  all $\nu\in \N_{\ge m}$  there exists a family of sections $(S_{\nu,j})_{j=1}^{M_{\nu}}\subset H^0(X,{\cal L}_{\nu})$ %$M_{\nu}:=N_{q_{\nu}}$
%$$(S_{\nu,j})_{j=1}^{M_{\nu}}\subset H^0(X,{\cal L}_{\nu})\,,\quad M_{\nu}:=N_{q_{\nu}}\,,$$
such that 
$
S_{\nu,j {|Z}}=\sigma_{\nu,j}\otimes(d\zeta)^{\nu}\,.
$
%for all $j=1,...,M_{\nu}$.
%with $\zeta\in H^0(X, {\cal O}(Z))$ such that $\div \zeta= Z$.
\end{Claim}
$Proof$. The proof of this claim goes by induction. The statement is obvious for $\nu=m$ by the assumption (A2). So we assume it true for $\nu$ and we prove it for $\nu+1$. 
We have 
\begin{eqnarray*}%\label{Ind-plurI}
{\cal L}_{\nu+1}
=
k_{\nu}F_m+(q_{\nu}+1)(K_X+Z)+\sum_{j=0}^{q_{\nu}+1}L_j
=
K_X+Z+{\cal L}_{\nu}+L_{q_{\nu}+1}\,,\nonumber
\end{eqnarray*}
if $q_{\nu}\le m-2$ and 
\begin{eqnarray*}%\label{Ind-plurII}
{\cal L}_{\nu+1}
\;=\;
(k_{\nu}+1)F_m+L_0
\;=\;
K_X+Z+{\cal L}_{\nu}+L_m\,,
\end{eqnarray*}
if $q_{\nu}=m-1$. So in all cases hold the induction formula 
$$
{\cal L}_{\nu+1}=K_X+Z+{\cal L}_{\nu}+L_{q_{\nu}+1}\,.%\qquad H_{\nu+1}=\Omega_X^{-1}\otimes h_Z\otimes H_{\nu}\otimes h_{L_{q_{\nu}+1}}\,.
$$ 
We will equip the line bundle ${\cal L}_{\nu}+L_{q_{\nu}+1}$ with an adequate singular hermitian metric with strictly positive curvature. For this purpose let $(\varepsilon_{\nu})_{\nu}, (\delta_{\nu})_{\nu}\subset (0,1)$ and consider the Hacon-McKernan decomposition \cite{Ha-Mc1}, \cite{Ha-Mc2}, \cite{Be-Pa}, \cite{Pa2}
\begin{eqnarray}\label{H-M-dec}
{\cal L}_{\nu}+L_{q_{\nu}+1}
&=&
(1-\varepsilon_{\nu}) L_{q_{\nu}+1}\,+\,\varepsilon_{\nu} L_{q_{\nu}+1}\nonumber
\\\nonumber
\\
&+&
(1-\delta_{\nu}){\cal L}_{\nu}\,+\,\delta_{\nu}(k_{\nu}F_m\,+\,{\cal L}_{q_{\nu}})\,,
\end{eqnarray}
in the case $q_{\nu}\le m-2$. The case $q_{\nu}=m-1$ will not present any difficulty.
Let $\tau_{\nu}\in (0,k_{\nu}^{-1})$.
According to Demailly's regularising process \cite{Dem1}, we can replace the current $\Theta$ with a family of closed and real $(1,1)$-currents with analytic singularities $\Theta_{\nu}\in \{\Theta\}$, $\Theta_{\nu}\ge -\tau_{\nu}\,\omega$, such that the restrictions $\Theta_{\nu |Z}$ are also well defined
and 
$$\lambda(\Theta_{\nu |Z})_z\;\le\; \lambda(\Theta_{ |Z})_z\,,
$$ 
for all $z\in Z$ and $\nu$. This combined with the condition \eqref{Sec-div} implies
\begin{eqnarray}\label{Sec-div-rg}
\div u- m\sum_{t=1}^{N'} \lambda(\Theta_{\nu|Z}, V_t)\, V_t\ge 0\,.
%\div u- \sum_{t=1}^{N'} m  \mu_t V_{Z,t}\ge 0\,,\qquad \mu_t:=\max\left\{\lambda_{Z,t}\,,\,\lambda(\Theta_{|Z}, V_{Z,t})/(2\pi) \right\}\,,
\end{eqnarray}
%still hold true. %Moreover the restirctions $\Theta_
Let now $\psi_{\nu}$ be a quasi-plurisubharmonic function with analytic singularities such that 
$$
2\pi\,m\,\Theta_{\nu}=i\,{\cal C}_{h_m}(F_m)+i\partial\bar\partial \psi_{\nu}\,,
$$
let $B_{\nu}:=\sum_{j=1}^{M_{\nu}}|S_{\nu,j}|_{H_{\nu}}^2$, let $\Phi_{\nu}:=\log B_{\nu}$ and set 
$$
\Psi_{\nu}\;:=\;
\left  \{
\begin{array}{lr}
(1-\varepsilon_{\nu})\varphi_{q_{\nu}+1}+(1-\delta_{\nu})\Phi_{\nu}+\delta_{\nu}\, k_{\nu} \,\psi_{\nu}\,,
\;\;\quad\mbox{if}\quad 
q_{\nu}\le m-2\,,
\\
\\
\varphi_m+\Phi_{\nu}\,,
\qquad\qquad\qquad\qquad\qquad\qquad\qquad\mbox{if}\quad 
q_{\nu}=m-1\,.
\end{array}
\right.
$$
We show now that for adequate choices of the parameters $\varepsilon_{\nu}, \delta_{\nu}$ the singular hermitian line bundle 
\begin{eqnarray}\label{Big-cnd}
({\cal L}_{\nu}+L_{q_{\nu}+1} \,,\,H_{\nu}\otimes h_{L_{q_{\nu}+1}}e^{-\Psi_{\nu}})\,,
\end{eqnarray}
is big and 
\begin{eqnarray}\label{L2-cnd}
\sigma_{\nu+1, j}\in H^0\left(Z\,,\,{\cal O}_Z({\cal L}_{\nu+1 {|Z}})\otimes{\cal I}(\Psi_{\nu {|Z}})\right)\,.
\end{eqnarray}
Then the conclusion of the claim \ref{Pau-lm} will follow by applying the claim \ref{OT-Bas-ext} to the section $\sigma_{\nu+1, j}\otimes (d\zeta)^{\nu}$ in order to obtain the required extensions $S_{\nu+1,j}$.
We distinguish again two cases.
\\
\\
{\bf Case $q_{\nu}\le m-2$}. By \eqref{H-M-dec} we infer the decomposition of the $(1,1)$-current	
\begin{eqnarray*}
&&i\,{\cal C}_{H_{\nu}}({\cal L}_{\nu})+i\,{\cal C}_{h_{L_{q_{\nu}+1}}}(L_{q_{\nu}+1})+i\partial\bar\partial \Psi_{\nu}
\\
\\
&=&(1-\varepsilon_{\nu})\left[i\,{\cal C}_{h_{L_{q_{\nu}+1}}}(L_{q_{\nu}+1})+i\partial\bar\partial \varphi_{q_{\nu}+1}\right]+\varepsilon_{\nu}\,i\,{\cal C}_{h_{L_{q_{\nu}+1}}}(L_{q_{\nu}+1})
\\
\\
&+&(1-\delta_{\nu})\Big[i\,{\cal C}_{H_{\nu}}({\cal L}_{\nu})+i\partial\bar\partial \Phi_{\nu}\Big]
+
\delta_{\nu}\Big[i\,{\cal C}_{H_{\nu}}({\cal L}_{\nu})
+
k_{\nu} i\partial\bar\partial\psi_{\nu}\Big]
\\
\\
&\ge&
\varepsilon_{\nu}\,i\,{\cal C}_{h_{L_{q_{\nu}+1}}}(L_{q_{\nu}+1})
+
\delta_{\nu}
\Big[2\pi\,k_{\nu}m\,\Theta_{\nu}+i\,{\cal C}_{H_{q_{\nu}}}({\cal L}_{q_{\nu}})\Big]
\\
\\
&\ge&
\varepsilon_{\nu}\,i\,{\cal C}_{h_{L_{q_{\nu}+1}}}(L_{q_{\nu}+1})+2\pi\,m\,\delta_{\nu}(1-k_{\nu} \tau_{\nu})\omega\,,
\end{eqnarray*}
by the assumption (A3). We infer that if $C_{\omega}>0$ is a constant such that $i\,{\cal C}_{h_j}(L_j)\ge -2\pi C_{\omega}\,\omega$ for all $j=1,...,m-1$ then the bundle \eqref{Big-cnd} is big as soon as 
\begin{eqnarray}\label{epsil-cond}
\varepsilon_{\nu}<m(1-k_{\nu} \tau_{\nu})\delta_{\nu}/C_{\omega}\,.
\end{eqnarray}
On the other hand the relation
$
\sigma_{\nu+1, j}=u^{k_{\nu}}\otimes s_{q_{\nu}+1,j}
$,
$j=1,...,M_{\nu+1}$ combined with the fact that
$$
B_{\nu\,{|Y}}=|d\zeta|^{2\nu}_{\omega, h_Z}\sum_{t=1}^{M_{\nu}}|u^{k_{\nu}}\otimes s_{q_{\nu},t}|_{\nu}^2\,,
$$
and with the definition of the constant $C$ implies
\begin{eqnarray}\label{Key-IntI}
I_{\nu+1}:=\int_Z|\sigma_{\nu+1, j}|^2_{\nu+1}\,e^{-\Psi_{\nu}}\le C'\int_Z|u|_{h'_m}^{2k_{\nu}\delta_{\nu}}e^{-(1-\varepsilon_{\nu})\varphi_{q_{\nu}+1}\,-\,\delta_{\nu} k_{\nu}\psi_{\nu}}\,,
\end{eqnarray}
with $h'_m:=\Omega_Z^{-m}\otimes h_{L_1}\otimes\cdots \otimes h_{L_m}$. We consider now the decomposition 
\begin{eqnarray}\label{decomp}
2\pi \,m\Theta_{\nu {|Z}}=2\pi[W_{\nu}]+\alpha_{\nu}+i\partial\bar\partial g_{\nu}\,,
\end{eqnarray}
with
$$
W_{\nu}:=m\sum_{t=1}^{N'}   \lambda(\Theta_{\nu |Z}, V_t)\, V_t\,,
$$
with $\alpha_{\nu}$ a smooth closed and real $(1,1)$-form and with $g_{\nu}$ a 
quasi-plurisub
harmonic function with analytic singularities such that $\lambda(g_{\nu}, V_t)=0$ for all $t=1,...,N'$. 
\\
In particular $g_{\nu}$ is not identically $-\infty$ over the sets $V_t$.
\\
Then the the decomposition \eqref{decomp} combined with the Lelong-Poincar\'{e} formula implies 
\begin{eqnarray*}
2\pi([\div u]-[W_{\nu}])
&=&
2\pi([\div u]-m\Theta_{\nu {|Z}})+\alpha_{\nu}+i\partial\bar\partial g_{\nu}
\\
\\
&=&\beta_{\nu}+i\partial\bar\partial f_{\nu}\,,
\end{eqnarray*}
with $\beta_{\nu}$ a smooth closed and real $(1,1)$-form and with
$$
f_{\nu}:=\log |u|_{h'_m}^2-\psi_{\nu}+g_{\nu}\,.
$$
The condition \eqref{Sec-div-rg} rewrites as $0\le \div u-W_{\nu}$.
We infer that $f_{\nu}$ is a quasi-plurisubharmonic function, thus bounded from above. We infer by \eqref{Key-IntI} the inequality 
\begin{eqnarray}\label{Key-IntII}
I_{\nu+1}\le C'\int_Ze^{-(1-\varepsilon_{\nu})\varphi_{q_{\nu}+1}\,+\,\delta_{\nu} k_{\nu}(f_{\nu}- g_{\nu})}\le C''\int_Ze^{-(1-\varepsilon_{\nu})\varphi_{q_{\nu}+1}\,-\,\delta_{\nu} k_{\nu}  g_{\nu}}\,.
\end{eqnarray}
On the other hand  
$$
\Lambda_{\nu}\;:=\;\sup_{z\in Z}\lambda(g_{\nu})_z\;<\;+\infty\,,
$$
since $Z$ is compact.
Thus the last integral in \eqref{Key-IntII} is convergent
for all values $\varepsilon_{\nu}\in (0,1)$ and 
\begin{eqnarray}\label{delt-cond}
0<\delta_{\nu} < 2(k_{\nu} \Lambda_{\nu})^{-1}\,,  
\end{eqnarray}
by the claim \ref{cv-int-Lel} and so the condition \eqref{L2-cnd} is satisfied in the case $q_{\nu}\le m-2$.
\\
\\
{\bf Case $q_{\nu}= m-1$}. In this case the condition \eqref{Big-cnd} is obviously satisfied. On the other hand the relation $\sigma_{\nu+1, j}=u^{k_{\nu}+1}\otimes s_{0,j}$ combined with the fact that
$$
B_{\nu\,{|Y}}=|d\zeta|^{2\nu}_{\omega, h_Z}\sum_{t=1}^{N_{m-1}}|u^{k_{\nu}}\otimes s_{m-1,t}|_{\nu}^2\,,
$$
and the definition of the constant $C$ implies 
\begin{eqnarray}\label{qishazot-intI}
I_{\nu+1}\le C'\int_Z|u|_{h'_m}^2e^{-\varphi_m}<+\infty\,,
\end{eqnarray}
The convergence follows from the condition \eqref{Sec-div} and the  fact that $\rho$ has bounded local potentials along $Z$. This concludes the proof of the claim \ref{Pau-lm}.\hfill$\Box$
\\
\\
{\bf End of the proof.} The claim \ref{Pau-lm} implies that  the singular hermitian line bundle 
$$
\left({\cal L}_{km}\,,\,H_{km}B_{km}^{-1}\right)\equiv \left(kF_m+A\,,\, h_m^k\otimes h_{A}\,B_{km}^{-1}\right)\,,
$$
is pseudoeffective.
So we have obtain the following;
\begin{eqnarray}
&&
i\,{\cal C}_{h_m}(F_m)\,+\,\frac{1}{k}\, i\,\partial\bar\partial \Phi_{km} \;\ge\; -\, \frac{1}{k}\,\omega_A\,,\label{fin-posit}
\\\nonumber
\\
&&
\frac{1}{k}\,\Phi_{km\,{|Z}}\;=\;\log |u|^2_{h'_m}\,+\,\frac{1}{k}\,\log\,\left(|d\zeta|^{2km}_{\omega, h_Z}\sum_{j=0}^{N_0}|s_{0,j}|^2_{h_A}\right)\,.\label{res-fin-pot}
\end{eqnarray}
Let $h_{mF}:=h_{L_1}\otimes\cdots \otimes  h_{L_m}$, let $\varphi_{\Delta}:=\frac{1}{m}\sum_{j=1}^{m-1}\varphi_j$ and set
$$
\Xi_k:=\frac{m-1}{mk}\,\Phi_{km}+\varphi_{\Delta}+\frac{1}{m}\varphi_m\,.
$$
Then the $\Q$-decomposition
$$
(m-1)(K_X+Z)+mL=\frac{m-1}{m}\,F_m+\Delta+\frac{1}{m}\,L_m\,,
$$
combined with the inequality \eqref{fin-posit} shows that the singular hermitian line bundle 
$$
\left((m-1)(K_X+Z)+mL \,,\, \Omega_X^{-(m-1)}\otimes h_Z^{m-1}\otimes h_{mL} \,e^{-\Xi_k}\right)\,,
$$
is big as soon as 
\begin{eqnarray}\label{k-cond}
k>(m-1)C_A/\varepsilon\,, 
\end{eqnarray}
with $\varepsilon\,,C_A\in \R_{>0}$ such that $\rho\ge \varepsilon \,\omega$ and $\omega_A\le 2\pi\, m C_A\,\omega$. On the other hand the expression \eqref{res-fin-pot} and the condition \eqref{Sec-div} imply
\begin{eqnarray*}
\int_Z|u|_{h'_m}^2e^{-\Xi_k} 
\;\le\;
C_k\int_Z|u|_{h'_m}^{2/m}e^{-\varphi_{\Delta}-\varphi_m/m}
\;\le\;
C'_k\int_Ze^{-\varphi_{\Delta}}\;<\;+\infty\,,
\end{eqnarray*}
since $h_{L_1}\otimes\cdots \otimes  h_{L_{m-1}}e^{-m\varphi_{\Delta}}$ is the canonical metric associated to the integral divisor $m\Delta$ and $\lambda_j<1$. % by the choice of the integer $m$.
In conclusion we can apply the claim \ref{OT-Bas-ext} to the section 
$$
u\otimes (d\zeta)^{m-1}\in H^0\Big(Z\,,\,K_Z+(m-1)(K_X+Z)+mL\Big)\,,
$$ 
in order to obtain the required lifting $U$ 
of the section $u$.\hfill$\Box$ 
%%%%%%%%%%%%%%%%%%%%%%%%%%%%%%%%%%%%%%%%%%%%%%%%%%%%%%%%%%%%%%%%%%%%%%%%%%%%%%%%%%%%%%%%%%%%%%%%%%%%%%%%%%%%%%%%%%%%%%%%%%%%%%%%%%%%%%%%%%%%%%%%%%%%%%%%%%%%%%%%%%%%%%%%%%%%
%%%%%%%%%%%%%%%%%%%%%%%%%%%%%%%%%%%%%%%%%%%%%%%%%%%%%%%%%%%%%%%%%%%%%%%%%%%%%%%%%%%%%%%%%%%%%%%%%%%%%%%%%%%%%%%%%%%%%%%%%%%%%%%%%%%%%%%%%%%%%%%%%%%%%%%%%%%%%%%%%%%%%%%%%%%%
%%%%%%%%%%%%%%%%%%%%%%%%%%%%%%%%%%%%%%%%%%%%%%%%%%%%%%%%%%%%%%%%%%%%%%%%%%%%%%%%%%%%%%%%%%%%%%%%%%%%%%%%%%%%%%%%%%%%%%%%%%%%%%%%%%%%%%%%%%%%%%%%%%%%%%%%%%%%%%%%%%%%%%%%%%%%
%%%%%%%%%%%%%%%%%%%%%%%%%%%%%%%%%%%%%%%%%%%%%%%%%%%%%%%%%%%%%%%%%%%%%%%%%%%%%%%%%%%%%%%%%%%%%%%%%%%%%%%%%%%%%%%%%%%%%%%%%%%%%%%%%%%%%%%%%%%%%%%%%%%%%%%%%%%%%%%%%%%%%%%%%%%%
%%%%%%%%%%%%%%%%%%%%%%%%%%%%%%%%%%%%%%%%%%%%%%%%%%%%%%%%%%%%%%%%%%%%%%%%%%%%%%%%%%%%%%%%%%%%%%%%%%%%%%%%%%%%%%%%%%%%%%%%%%%%%%%%%%%%%%%%%%%%%%%%%%%%%%%%%%%%%%%%%%%%%%%%%%%%
%%%%%%%%%%%%%%%%%%%%%%%%%%%%%%%%%%%%%%%%%%%%%%%%%%%%%%%%%%%%%%%%%%%%%%%%%%%%%%%%%%%%%%%%%%%%%%%%%%%%%%%%%%%%%%%%%%%%%%%%%%%%%%%%%%%%%%%%%%%%%%%%%%%%%%%%%%%%%%%%%%%%%%%%%%%%
\newpage
\subsection{A perturbed extension statement}
The Ein-Popa extension result \cite{Ei-Po} previously explained modifies quite directly in a perturbed extension statement due to Paun \cite{Pa2}. We explain now this statement.
For any $\Q$-line bundle/divisor $E$ we fix a smooth form $\theta_E\in c_1(E)$. %associated to the smooth hermitian line bundles $(E,h_E)$ considered in the proof of the lemma \ref{Hk-Mk-lem}.
%In the following statement we tacitly assume that the line bundles ${\cal O}_X(Z)$, ${\cal O}_X(Z_j)$ are equiped with smooth hermitian metrics.
We observe the following quite elementary fact.
%%%%%%%%%%%%%%%%%%%%%%%%%%%%%%%%%%%%%%%%%%%%%%%%%%%%%%%%%%%%%%%%%%%%%%%%%%%%%%%%%%%%%%%%%%%%%%%%%%%%%%%%%%%%%%%%%%%%%%%%%%%%%%%%%%%%%%%%%%%%%%%%%%%%%%%%%%%%%%%%%%%%%%%%%%%%
%%%%%%%%%%%%%%%%%%%%%%%%%%%%%%%%%%%%%%%%%%%%%%%%%%%%%%%%%%%%%%%%%%%%%%%%%%%%%%%%%%%%%%%%%%%%%%%%%%%%%%%%%%%%%%%%%%%%%%%%%%%%%%%%%%%%%%%%%%%%%%%%%%%%%%%%%%%%%%%%%%%%%%%%%%%%
%%%%%%%%%%%%%%%%%%%%%%%%%%%%%%%%%%%%%%%%%%%%%%%%%%%%%%%%%%%%%%%%%%%%%%%%%%%%%%%%%%%%%%%%%%%%%%%%%%%%%%%%%%%%%%%%%%%%%%%%%%%%%%%%%%%%%%%%%%%%%%%%%%%%%%%%%%%%%%%%%%%%%%%%%%%%
\begin{Claim}\label{Univ-Amp}
Let $A_0$ be an ample line bundle over a complex projective variety $X$ of complex dimension $n$, let $\omega\in c_1(A_0)$, $\omega>0$, let  $Z$, $Z_j\subset X$, $j=1,...,N$ be irreducible divisors and let $D$ be a $\Q$-divisor over $X$. Let also $C_0\in \N_{>0}$ such that
$$
\theta_Z\,,\theta_{Z_j}\,,\theta_D\,,\theta_{K_X}\,,\frac{n-1}{\pi}\,i\partial\bar\partial \log\dist_{\omega}(x,\cdot)\;\ge\;-C_0\,\omega\,,
$$
for all $j=1,...,N$ and $x\in X$.
Then for any holomorphic $\Q$-line bundle $R$ as in the statement of the lemma \ref{Hk-Mk-lem}, any $m\in \N_{>1}$ such that $mD$, $mR$ are integral and any subset 
$$
{\cal S}\subset \{Z_j\,:\,j=1,...,N\}\times \{1,...,m-1\}\,,
$$ 
the family of holomorphic line bundles $(L_k)_{k=1}^m$ defined by 
$$
L_k:={\cal O}_X(\Delta_k)\,,\qquad \Delta_k:=\sum_{Z\in {\cal S}_k}Z\,,\qquad {\cal S}_k:=\{Z_j\,:\,j=1,...,N\}\times \{k\}\,,
$$
for all $k=1,...,m-1$ and $L_m:={\cal O}_X(mD)+m\,R$, satisfies the properties $(AI)$, $I=1,2,3$ in the proof of the lemma \ref{Hk-Mk-lem} with respect to 
$$
A:=m[2+(N+3)C_0]A_0\,.
$$
\end{Claim}
$Proof$.  The inequality $\theta_{L_k}\ge -NC_0\,\omega$ for all $k=1,...,m-1$ implies
$$
\theta_{{\cal L}_q}\;\ge\;\theta_A-(m-1)(N+2)C_0\,\omega\,,\quad \forall q=0,...,m-1\,.
$$
For $q=m$ hold the inequality 
$$
\Theta_{{\cal L}_m}\;\ge\;\theta_A-(m-1)(N+2)C_0\,\omega-mC_0\omega\,,
$$
where $\Theta_{{\cal L}_m}\in c_1({\cal L}_m)$ is a current with bounded potentials along $Z$.
On the other hand the Kawamata-Viehweg-Nadel vanishing theorem and the claim \ref{OT-Bas-ext} imply that the properties $(A1)$ and $(A2)$ in the proof of the lemma \ref{Hk-Mk-lem} are satisfied with respect to $A$ in the statement of the claim \ref{Univ-Amp}. This choice of $A$ satisfies also the property $(A3)$.\hfill$\Box$
%%%%%%%%%%%%%%%%%%%%%%%%%%%%%%%%%%%%%%%%%%%%%%%%%%%%%%%%%%%%%%%%%%%%%%%%%%%%%%%%%%%%%%%%%%%%%%%%%%%%%%%%%%%%%%%%%%%%%%%%%%%%%%%%%%%%%%%%%%%%%%%%%%%%%%%%%%%%%%%%%%%%%%%%%%%%
%%%%%%%%%%%%%%%%%%%%%%%%%%%%%%%%%%%%%%%%%%%%%%%%%%%%%%%%%%%%%%%%%%%%%%%%%%%%%%%%%%%%%%%%%%%%%%%%%%%%%%%%%%%%%%%%%%%%%%%%%%%%%%%%%%%%%%%%%%%%%%%%%%%%%%%%%%%%%%%%%%%%%%%%%%%%
%%%%%%%%%%%%%%%%%%%%%%%%%%%%%%%%%%%%%%%%%%%%%%%%%%%%%%%%%%%%%%%%%%%%%%%%%%%%%%%%%%%%%%%%%%%%%%%%%%%%%%%%%%%%%%%%%%%%%%%%%%%%%%%%%%%%%%%%%%%%%%%%%%%%%%%%%%%%%%%%%%%%%%%%%%%%%%%%%%%%%%%%%%%%%%%%%%%%%%%%%%%%%%%%%%%%%%%%%%%%%%%%%%%%%%%%%%%%%%%%%%%%%%%%%%%%%%%%%%%%%%%%%%%%%%%%%%%%%%%%%%%%%%%%%%%%%%%%%%%%%%
%This implies the following analytic version of the  Hacon-McKernan and Ein-Popa  extension result for pluri-twisted holomorphic sections of singular hermitian line 
%%%%%%%%%%%%%%%%%%%%%%%%%%%%%%%%%%%%%%%%%%%%%%%%%%%%%%%%%%%%%%%%%%%%%%%%%%%%%%%%%%%%%%%%%%%%%%%%%%%%%%%%%%%%%%%%%%%%%%%%%%%%%%%%%%
%%%%%%%%%%%%%%%%%%%%%%%%%%%%%%%%%%%%%%%%%%%%%%%%%%%%%%%%%%%%%%%%%%%%%%%%%%%%%%%%%%%%%%%%%%%%%%%%%%%%%%%%%%%%%%%%%%%%%%%%%%%%%%%%%%
\newpage
\begin{corol}\label{Pert-Hk-Mk-lem}
Let $X$ be a complex projective manifold,
let $Z\subset X$ be a smooth irreducible hypersurface,
let $A_0$ be an ample line bundle over $X$, 
let $\omega\in c_1(A_0)$, $\omega>0$  
and let $L$ be a  holomorphic $\Q$-line bundle over $X$ which admits a decomposition as 
$$
L={\cal O}_X(\Delta)+{\cal O}_X(D)+R\,,
$$
where;
\\
$\blacktriangleright$ 
$
\Delta=\sum_{j=1}^N \lambda_j\, Z_j
$
is a divisor over $X$ with $\lambda_j\in \Q\cap [0,1)$ and $Z_j\subset X$ distinct irreducible smooth hypersurfaces with normal crossing intersection with $Z$ such that $Z\cap Z_j\cap Z_l =\emptyset$ for all $j\not = l$,
\\
\\
$\blacktriangleright$
$D$ is an effective $\Q$-divisor over $X$ such that $Z$ is not one of its components and the components $(\Gamma_p)_{p=1}^Q$ of the restricted  divisor $D_{|Z}$ does not intersect the irreducible components $(V_t)_{t=1}^{N'}$ of the family $(Z_j\cap Z)_{j=1}^N$,
\\
\\
$\blacktriangleright$
$R$ is a holomorphic $\Q$-line bundle over $X$ such that 
there exists a K\"{a}hler current
$\rho\in c_1(R)$ with $\rho\ge \varepsilon\omega$, $\varepsilon\in \R_{>0}$ and with bounded local potentials along $Z$. 
%well defined restricton $\theta_{|Z}$.% such that ${\cal I}(m^{-1}\theta_{|Z})={\cal O}_Z$,
\\
\\
$\bullet$ Let  $C_0\in \N_{>0}$ as in the statement of the claim \ref{Univ-Amp}, let 
\begin{eqnarray*}
C_1:=2+(N+3)C_0\,,\qquad
C_2:=NC_0C_1\,,\qquad
\lambda&:=&\max_{1\le j\le N} \lambda_j\,.
\end{eqnarray*}
$\bullet$ Let $m\in\N_{>1}$,  such that $m\Delta$, $mD$ are integral, $mR$ is a holomorphic line bundle and
$$
m\;\ge\;\frac{1}{2C_2(1-\lambda)\lceil1/\varepsilon\rceil}\,.
$$
$\bullet$ Let $V:=\sum_{t=1}^{N'}V_t$, let $\Gamma:=\sum_{p=1}^Q\Gamma_p$ and let $\eta\in \R_{>0}$ such that $\eta<1/\mult(\Gamma)$.
\\
\\
Assume the existence of a closed $(1,1)$-current $\Theta \in c_1(K_X+Z+L)$
%$$\Theta \in c_1(K_X+Z+L)\,,$$ 
with analytic singularities and with well defined restriction  $\Theta_{|Z}$ such that 
\begin{eqnarray}\label{err-pos-curr}
\Theta\;\ge\; -\;\frac{1}{2C_1\lceil1/\varepsilon\rceil}\;\frac{1}{m}\;\omega\,.
\end{eqnarray}
Then  for any 
$
u\in H^0\Big(Z\,,\,m\left(K_Z+ L_{{|Z}}\right)\Big)%\,,
$
with the vanishing property
\begin{eqnarray}\label{prt-Sec-div}
\div u\;-\; m\left(\,\sum_{t=1}^{N'} \lambda(\Theta_{|Z}\,,V_t)\, V_t\;+\;D_{|Z}\right)\;\ge\; -\;\frac{1}{3C_2\lceil1/\varepsilon\rceil}\,V\;-\;\eta\Gamma	\,,\quad
%\div u- \sum_{t=1}^{N'} m  \mu_t V_{Z,t}\ge 0\,,\qquad \mu_t:=\max\left\{\lambda_{Z,t}\,,\,\lambda(\Theta_{|Z}, V_{Z,t})/(2\pi) \right\}\,,
\end{eqnarray}
there exists a section 
$$
U\in H^0\Big(X\,,\,m(K_X+Z+L)\Big)\,,\qquad
U_{{|Z}}=u\otimes(d\zeta)^m\,,
$$
with $\zeta\in H^0(X, {\cal O}(Z))$ such that $\div \zeta= Z$.
\end{corol}
$Proof$.
We repeat the proof of the lemma \ref{Hk-Mk-lem} with some very little modifications. 
The data $(\lambda_j)_j$ determines a set $\cal S$ as in the statement of the claim \ref{Univ-Amp}. 
Thus the conditions $(AI)$, $I=1,2,3$ in the proof of 
the lemma \ref{Hk-Mk-lem} are satisfied with respect to $A$ in statement of the claim \ref{Univ-Amp}. 
We perform the induction of the claim \ref{Pau-lm} for the steps $\nu=m,...,\bar km$, with $\bar k:=m\,C_1\lceil1/\varepsilon\rceil$.
\\
\\
{\bf The case $q_{\nu}\le m-2$.}
We replace the currents $\Theta_{\nu}$ in the proof of the lemma \ref{Hk-Mk-lem} with 
the current $\Theta\ge -\tau\,\omega$, 
$$
\tau:=\frac{1}{2\bar k}\,,
$$ 
and we reconsider the conditions needed 
for the parameters $\varepsilon_{\nu}\equiv\bar\varepsilon>0$, $\delta_{\nu}\equiv\bar\delta>0$. 
With the notations of the claim \ref{Univ-Amp} hold the inequality
$
\theta_{L_k}\ge-NC_0\,\omega
$.
We infer that in our setting the condition \eqref{epsil-cond} on the bigness of the singular hermitian line bundle \eqref{Big-cnd} becomes
$$
0\,<\,\bar\varepsilon\,<\,\frac{m(1-k_{\nu}\tau)\bar\delta}{NC_0}\,.
$$
We observe that the inequality $1-k_{\nu}\tau>0$ is satisfied for all $\nu=m,...,\bar k m$ by our definition of  $\tau>0$. 
So a first condition on $\bar\varepsilon$ is
$$
\bar\varepsilon\,<\,\frac{m(1-\bar k\tau)\bar\delta}{NC_0}\,.
$$
Let now $\psi$, $W$, $\alpha$ and $g$ correspond respectively to $\psi_{\nu}$, $W_{\nu}$, $\alpha_{\nu}$ and $g_{\nu}$ in the proof of the claim \ref{Pau-lm} and let $\varphi_V$, $\varphi_{\Gamma}$ such that 
$$
2\pi[V]=\theta_V+i\partial\bar\partial \varphi_V\,,\qquad 2\pi[\,\Gamma\,]=\theta_{\Gamma}+i\partial\bar\partial \varphi_{\Gamma}\,,
$$
for some smooth $(1,1)$-forms $\theta_V$ and $\theta_{\Gamma}$. 
Let 
$$
\mu:=\frac{1}{3C_2\lceil1/\varepsilon\rceil}\,.
$$
This definition implies the inequality
\begin{eqnarray}\label{mu-cond}
\mu<m\min\left\{\frac{1/\bar k-\tau}{NC_0}\,,\,1-\lambda\right\}\,,
\end{eqnarray}
by our choice of $m$.
By the vanishing condition \eqref{prt-Sec-div} and the Lelong-Poincar\'{e} formula we infer
\begin{eqnarray*}
0&\le&
2\pi\left([\div u-W+\mu V+\eta\Gamma]\right) 
\\
\\
&=&
2\pi\left([\div u+\mu V+\eta\Gamma]-m\Theta_{|Z}\right)+\alpha+i\partial\bar\partial g
\\
\\
&=&\beta+i\partial\bar\partial f\,,
\end{eqnarray*}
with $\beta$ a smooth $(1,1)$-form and with
$$
f\;:=\;\log |u|_{h'_m}^2\,-\,\psi\,+\,g\,+\,\mu\varphi_V\,+\,\eta\varphi_{\Gamma}\,,
$$
quasi-plurisubharmonic, thus bounded from above. We infer
\begin{eqnarray}\label{Prt-Key-IntII}
I_{\nu+1}
&\le&
C'\int_Ze^{-(1-\bar\varepsilon)\varphi_{q_{\nu}+1}\,+\,\bar\delta k_{\nu}(f- g-\mu\varphi_V-\eta\varphi_{\Gamma})}\nonumber
\\\nonumber
\\
&\le&
C''\int_Ze^{-(1-\bar\varepsilon)\varphi_{q_{\nu}+1}\,-\,\bar\delta k_{\nu}\mu\varphi_V\,-\,\bar\delta k_{\nu}  (g+\eta\varphi_{\Gamma})}\,.
\end{eqnarray}
Let
$$
\Lambda_{\eta}:=\sup_{z\in Z}\lambda(g+\eta\varphi_{\Gamma})_z\,<\,+\infty\,.
$$
By the claim \ref{cv-int-Lel} the integral \eqref{Prt-Key-IntII} is finite if $\bar\delta k_{\nu}\mu<\bar \varepsilon$ and $\bar\delta\,<\,2(k_{\nu}\Lambda_{\eta})^{-1}$
for all 
$\nu=m,...,\bar k m$. So we take $\bar \varepsilon$ and $\bar\delta$ such that 
$$
\bar\delta \bar k\mu\,<\,\bar \varepsilon  \,<\, \frac{m(1-\bar k\tau)\bar\delta}{NC_0}\,,\qquad            0\,<\,\bar\delta\,<\,2(\bar k\Lambda_{\eta})^{-1}\,.
$$
The existence of $\bar \varepsilon$ follows from the inequality \eqref{mu-cond}.
\\
\\
{\bf The case $q_{\nu}= m-1$.} We consider the decomposition $\varphi_m=\varphi_{mD}+\varphi_{m\rho}$, where $\varphi_{mD}$ and $\varphi_{m\rho}$ are potentials corresponding respectively to the closed positive currents $2\pi[mD]$ and $2\pi m\rho$. The vanishing condition \eqref{prt-Sec-div} and the Lelong-Poincar\'{e} formula imply
$$
0\;\le\;
2\pi\left([\div u-mD_{{|Z}}+\mu V+\eta\Gamma]\right)\;=\;\tilde\beta+ i\partial\bar\partial \tilde f\,,
$$
with $\tilde\beta$ a smooth $(1,1)$-form and with
$$
\tilde f:=\log |u|_{h'_m}^2-\varphi_{mD}+\mu\,\varphi_V+\eta\,\varphi_{\Gamma}\,,
$$
quasi-plurisubharmonic, thus bounded from above. We infer
$$
I_{\nu+1}\le \tilde C' \int_Z|u|_{h'_m}^2e^{-\varphi_m}
=
C' \int_Ze^{\tilde f -\varphi_{m\rho}-\mu\,\varphi_V-\eta\,\varphi_{\Gamma}}\le \tilde C''\int_Z e^{-\mu\,\varphi_V-\eta\,\varphi_{\Gamma}} <+\infty\,,
$$
since $\mu<1$, since the singular part of $\varphi_V$ does not intersect with the singular part of $\varphi_{\Gamma}$ and since $\varphi_{m\rho}$ is bounded along $Z$ by assumption.
\\
\\
{\bf End of the proof.}
The constant $C_A>0$ in the proof of the lemma \ref{Hk-Mk-lem} corresponds to $C_1$. We infer that the condition \eqref{k-cond} becomes $k>(m-1)C_1/ \varepsilon$, which is satisfied by our choice of the integer $\bar k$. On the other hand
\begin{eqnarray*}
\int_Z|u|_{h'_m}^2e^{-\Xi_{\bar k}} 
&\le&	
C_k\int_Z|u|_{h'_m}^{2/m}e^{-\varphi_{\Delta}\,-\,\varphi_m/m}
\\
\\
&=&
C_k\int_Ze^{-\varphi_{\Delta}\,+\,(\tilde f \,-\,\varphi_{m\rho}\,-\,\mu\varphi_V\,-\,\eta\,\varphi_{\Gamma})/m}
\\
\\
&\le&
C'_k\int_Ze^{-\varphi_{\Delta}\,-\,\mu\varphi_V/m\,-\,\eta\,\varphi_{\Gamma}/m}\;<\;+\infty\,,
\end{eqnarray*}
since $\lambda_j+\mu/m<1$ for all $j=1,...,N$ 
by the inequality \eqref{mu-cond} and since the singular part of $\varphi_{\Delta}+\mu\varphi_V/m$ 
does not intersect with the singular part of $\varphi_{\Gamma}$ 
by our assumption on the components of the divisor $D_{{|Z}}$.
\hfill $\Box$
\newpage
%%%%%%%%%%%%%%%%%%%%%%%%%%%%%%%%%%%%%%%%%%%%%%%%%%%%%%%%%%%%%%%%%%%%%%%%%%%%%%%%%%%%%%%%%%%%%%%%%%%%%%%%%%%%%%%%%%%%%%%%%%%%%%%%%%
%%%%%%%%%%%%%%%%%%%%%%%%%%%%%%%%%%%%%%%%%%%%%%%%%%%%%%%%%%%%%%%%%%%%%%%%%%%%%%%%%%%%%%%%%%%%%%%%%%%%%%%%%%%%%%%%%%%%%%%%%%%%%%%%%%
%%%%%%%%%%%%%%%%%%%%%%%%%%%%%%%%%%%%%%%%%%%%%%%%%%%%%%%%%%%%%%%%%%%%%%%%%%%%%%%%%%%%%%%%%%%%%%%%%%%%%%%%%%%%%%%%%%%%%%%%%%%%%%%%%%
%%%%%%%%%%%%%%%%%%%%%%%%%%%%%%%%%%%%%%%%%%%%%%%%%%%%%%%%%%%%%%%%%%%%%%%%%%%%%%%%%%%%%%%%%%%%%%%%%%%%%%%%%%%%%%%%%%%%%%%%%%%%%%%%%%
%%%%%%%%%%%%%%%%%%%%%%%%%%%%%%%%%%%%%%%%%%%%%%%%%%%%%%%%%%%%%%%%%%%%%%%%%%%%%%%%%%%%%%%%%%%%%%%%%%%%%%%%%%%%%%%%%%%%%%%%%%%%%%%%%%
%%%%%%%%%%%%%%%%%%%%%%%%%%%%%%%%%%%%%%%%%%%%%%%%%%%%%%%%%%%%%%%%%%%%%%%%%%%%%%%%%%%%%%%%%%%%%%%%%%%%%%%%%%%%%%%%%%%%%%%%%%%%%%%%%%
%%%%%%%%%%%%%%%%%%%%%%%%%%%%%%%%%%%%%%%%%%%%%%%%%%%%%%%%%%%%%%%%%%%%%%%%%%%%%%%%%%%%%%%%%%%%%%%%%%%%%%%%%%%%%%%%%%%%%%%%%%%%%%%%%%
%%%%%%%%%%%%%%%%%%%%%%%%%%%%%%%%%%%%%%%%%%%%%%%%%%%%%%%%%%%%%%%%%%%%%%%%%%%%%%%%%%%%%%%%%%%%%%%%%%%%%%%%%%%%%%%%%%%%%%%%%%%%%%%%%%
%%%%%%%%%%%%%%%%%%%%%%%%%%%%%%%%%%%%%%%%%%%%%%%%%%%%%%%%%%%%%%%%%%%%%%%%%%%%%%%%%%%%%%%%%%%%%%%%%%%%%%%%%%%%%%%%%%%%%%%%%%%%%%%%%%
%%%%%%%%%%%%%%%%%%%%%%%%%%%%%%%%%%%%%%%%%%%%%%%%%%%%%%%%%%%%%%%%%%%%%%%%%%%%%%%%%%%%%%%%%%%%%%%%%%%%%%%%%%%%%%%%%%%%%%%%%%%%%%%%%%
\section{Applications to the Non-Vanishing}\label{sng-ext-rst}
In this section we will discuss the issues related with the application of the perturbed extension statement of Ein-Popa type to the following fundamental Non-Vanishing result due to the team BCHM \cite{BCHM}. %We explain here the setting.
Let $X$ be a complex projective manifold. We will consider the Neron-Severi lattice and the associated real Neron-Severi space%$\Q$ or $\R$-Neron-Severi spaces
\begin{eqnarray*}
\NS(X)&:=&H^{1,1}(X,\R)\cap (H^2(X,\Z)/\{\mbox{torsion}\})\,, 
%\\\\%\NS_{\Q}(X)&:=&\NS(X)\otimes _{\Z}\Q\simeq H^{1,1}(X,\R)\cap H^2(X,\Q)\,,
\\
\\
\NS_{\R}(X)&:=&\NS(X)\otimes _{\Z}\R\,.
\end{eqnarray*}
%With this notations we have the following result.
\begin{theorem}{\bf(Non-Vanishing result).} \label{A-NV}
\\
Let $X$ be a complex projective manifold and let $\alpha\in \NS_{\R}(X)$ admitting a K\"{a}hler current $\theta\in \alpha$ with ${\cal I}(\theta)={\cal O}_{X}$ and such that the $(1,1)$-cohomology class $c_1(K_X)+\alpha$ is pseudoeffective.
Then there exist a non zero effective $\R$-divisor $D$ such that
$$
[D]\in c_1(K_X)+\alpha\,.
$$
\end{theorem}
\subsection{Shokurov's modification of the manifold}
The following construction is essentially due to Shokurov. The slight modification explained here is similar to the one explained in \cite{Pa2}.
\begin{lem}\label{Sokur-const}
Under the assumptions of the Non-Vanishing theorem  \ref{A-NV} and under the assumption that the Non-Vanishing theorem \ref{A-NV} hold true over any complex projective manifold of dimension $\dim_{_{\C}}X-1$ there exist; 
\\
$\bullet$ a birational morphism $\mu:\widehat X\rightarrow X$ of complex projective manifolds and a K\"{a}hler form $\widehat \omega$ over $\widehat X$,
\\
$\bullet$ finitely many smooth and irreducible hypersurfaces $Z\,,Z_j\subset \widehat X$, $j\in I$ as in the in the statement of the lemma \ref{Hk-Mk-lem} such that  the sets $V_j:=Z_j\cap Z$ are irreducible for all $j\in I$. 
\\
$\bullet$ a $(1,1)$-cohomology class  $\widehat\alpha:=\{[\Delta]\}+\kappa$, with %$\Delta=\sum _{j\in I}a_j\,Z_j$
$$
%\widehat\alpha:=\{[\Delta]\}+\kappa\,,\qquad 
\Delta=\sum _{j\in I}a_j\,Z_j
\,,\qquad 
(a_j)_{j\in I}\subset (0,1)
$$ 
and with $\kappa \in \NS_{{\R}}(\widehat X)$ a K\"{a}hler class,
such that the following hold:
\\
$\blacktriangleright$ there exist a sequence $(\varepsilon_l)_l\subset (0,1)$, $\varepsilon_l\downarrow 0$ as $l\rightarrow +\infty$ and closed $(1,1)$-currents 
$$
\widehat\Theta_l\in c_1(K_{\widehat X}+Z)+\widehat{\alpha}\,,\qquad \widehat\Theta_l\;\ge\;-\varepsilon_l\,\widehat\omega\,,
$$ 
with well defined restriction $\widehat\Theta_{l |Z}$,
\\
$\blacktriangleright$ there exist an effective $\R$-divisor $G$ over $Z$ such that  $[G]\in c_1(K_Z)+\widehat\alpha_{ |Z}$ and such that for all $l\in \N$,
\begin{eqnarray}\label{NV-Vansh-cnd}
G\;-\;\sum_{j\in I_{-}}\lambda(\widehat\Theta_{l |Z},V_j)\,V_j\;-\; \sum _{j\in I_{+}}a_j\,V_j\;\ge \;0\,, 
\end{eqnarray}
with respect to a decomposition $I=I_{+}\amalg I_{-}$ independent of $l$,
\\
$\blacktriangleright$ there exist an effective $\R$-divisor  $F$ over $\widehat X$ and $r\in \R_{>1}$ such that 
\begin{eqnarray}\label{Zot-Shok}
c_1(K_{\widehat X}+Z)+\widehat{\alpha}+\{[F]\}\;=\;r\,\mu^*\Big(c_1(K_X)+\alpha\Big)+\{[E]\}\,,
\end{eqnarray}
with  ${\cal O}_{\widehat X}(E):=K_{\widehat X}-\mu^*K_X$.
\end{lem}
{\bf Proof.}
\\
{\bf (A) Set up and notations.} 
\\
Let $\omega>0$ be a K\"{a}hler form inside a rational class and let $\varepsilon_0\in \R_{>0}$ such that $\theta\ge \varepsilon_0\,\omega$.
By means of Demailly's regularisation process we can assume without lost of generality that the local potentials of the current $\theta$ are with regular analytic singularities and its Lelong numbers are rational. In fact  Demailly's regularisation process preserves the conditions $\theta\ge \varepsilon'\omega$, for some $\varepsilon'\in (0,\varepsilon_0)$ and  ${\cal I}(\theta)={\cal O}_{X}$. Let $\Theta\in c_1(K_X)+\alpha$ be a closed positive $(1,1)$-current and 
fix now a point $x_0\in X$ such that the local potentials of $\theta$ and $\Theta$ are bounded at the point $x_0$. By standard facts one infers the existence of $p\in \N_{>0}$ and of a closed positive current 
$$
H\in \beta:=p\Big(c_1(K_X)+\alpha\Big)+\alpha\,,
$$ 
with regular analytic singularities and rational Lelong numbers such that 
$$
\lambda(H,x_0)\ge n+1\,,
$$
$n=\dim_{\C}X$. By means of the Hironaka desingularisation result we can find a  birational morphism $\mu:\widehat X\rightarrow X$ of projective manifolds  which factors through the blow up map of $x_0$ and such that; 
\begin{eqnarray*}
K_{\widehat X} =\mu^{*}K_X+E\,,\qquad
E:=\sum_{j\in \bar J}e_j Z_j\,,\qquad
e_j\in \Z_{\ge 0}\,,\qquad
Z_{j_0}:=\mu^{-1}(x_0)\,,
\end{eqnarray*}
with $(Z_j)_{j\in \bar J}$ a finite family of distinct and smooth irreducible hypersurfaces with simple normal crossings such that $Z_j\cap Z_k\cap Z_l=\emptyset$ for all $j\not=k\not=l\not=j$, such that $Z_j\cap Z_k$ are irreducible for all $j,k$ and such that,
\begin{eqnarray*}
\mu^{*}\theta&=&\sum_{j\in \bar J}\theta_j[Z_j]+R_{\theta}\,,\quad
\theta_j\in \Q_{\ge 0}\,,\quad
\theta_{j_0}=0\,,\quad
R_{\theta}\ge \varepsilon'\mu^{*}\omega \quad \mbox{smooth,}
\\
\\
\mu^{*}H&=&\sum_{j\in \bar J}h_j [Z_j]+R_H\,,\quad
h_j\in \Q_{\ge 0}\,,\quad
h_{j_0}\ge n+1\,,\quad
R_H\ge 0
\quad \mbox{smooth.}
\end{eqnarray*}
Notice also that $e_{j_0}=n-1$ since $\mu$ factors through the blow up map of $x_0$. 
Consider now the decomposition %(allso used in Siu and Paun)
$$
\mu^{*}\Theta=\sum_{j\in \bar J}\lambda_j[Z_j]+R_{\Theta}\,,\quad 
\lambda_j:=\lambda(\mu^{*}\Theta\,,Z_j)\in \R_{\ge 0}
$$
whith $R_{\Theta}$ a closed positive $(1,1)$-current such that $\lambda(R_{\Theta}\,,Z_j)=0$ by construction. Notice that $\lambda_{j_0}=0$ by our choice of the point $x_0$. Let $\widehat \omega>0$ be a K\"{a}hler form over $\widehat X$. Let now $\tau\in \Q\cap (0,1)$ and apply the $\mu^*$-functor to the decomposition
$$
c_1(K_X)+\tau\beta+(1-\tau)\alpha\;=\;(\tau p+1)\Big(c_1(K_X)+\alpha\Big)\,,
$$
(also considered in Paun \cite{Pa2})
we infer the cohomology identity
\begin{eqnarray}\label{key-id}
c_1(K_{\widehat X})+\mu^*\Big(\tau\beta+(1-\tau)\alpha\Big)-\{[E]\}=(\tau p+1)\mu^*\Big(c_1(K_X)+\alpha\Big)\,.
\end{eqnarray}
Let $\zeta_j\in \{Z_j\}$ be a smooth representative, let $\delta_j\in \Q_{>0}$ and set
\begin{eqnarray*}
\widehat\omega_{\delta}&:=&\tau R_H+(1-\tau)R_{\theta}-\sum_{j\in \bar J}\delta_j\,\zeta_j\;\ge\; \varepsilon'\mu^*\omega-\sum_{j\in \bar J}\delta_j\,\zeta_j\,,
\\
\\
a_j&:=&\tau(h_j-\theta_j-p\lambda_j)-\lambda_j+\theta_j-e_j+\delta_j\,,
\\
\\
b_j&:=&\tau(h_j-\theta_j)+\theta_j-e_j+\delta_j\,,
\end{eqnarray*}
We remind that there exist $\delta_{\omega, \zeta}\in \R_{>0}$ such that 
$$
\varepsilon'\mu^*\omega-\sum_{j\in \bar J}\delta_j\,\zeta_j>0\,,
$$ 
for all $\delta_j\in (0,\delta_{\omega, \zeta})$. We infer by \eqref{key-id} the identities 
\begin{eqnarray}
c_1(K_{\widehat X})+ \sum_{j\in \bar J}a_j\,\{[Z_j]\}+\{\widehat\omega_{\delta}\}&=&(\tau p+1)\{R_{\Theta}\}\,,\label{key-idI}
\\\nonumber
\\
c_1(K_{\widehat X})+ \sum_{j\in \bar J}b_j\,\{[Z_j]\}+\{\widehat\omega_{\delta}\}&=&(\tau p+1)\mu^*\Big(c_1(K_X)+\alpha\Big)\,.\label{key-idII}
\end{eqnarray}
{\bf (B) Tie breaking.} 
\\
The assumption ${\cal I}(\theta)={\cal O}_X$ is equivalent to the inequality $\theta_j-e_j<1$ for all $j\in \bar J$, which implies
\begin{eqnarray}\label{Pau-Intg}
-\lambda_j+ \theta_j-e_j<1\,,\qquad \forall j\in \bar J\,.
\end{eqnarray}
On the other hand the inequality $h_{j_0}-\theta_{j_0}-p\lambda_{j_0}=h_{j_0}\ge n+1>0$ implies
$$
j_0\in I:=\{j\in \bar J\,:\,h_j-\theta_j-p\lambda_j>0\}\not =\emptyset\,.
$$
Let define now 
$$
c_j^{\delta}:=\frac{1+\lambda_j- \theta_j+e_j-\delta_j}{\,h_j-\theta_j-p\lambda_j}>0\,,
$$
for all $j\in I$. We choose $\delta_j\in \Q_{>0}$, $j\in \bar J$ such that 
\begin{eqnarray}\label{Tie-bk-1}
-\lambda_j+ \theta_j-e_j+\delta_j<1\,,\qquad \forall j\in \bar J\,, 
\end{eqnarray}
there exist a unique $j_1\in I$ such that 
\begin{eqnarray}\label{Tie-bk-2}
0\;<\;\tau:=\min_{j\in I} c_j^{\delta}\;<\;c_j^{\delta}\,,\qquad \forall j\in I\smallsetminus \{j_1\}
\end{eqnarray}
and such that
\begin{eqnarray}\label{Tie-bk-3}
\frac{\varepsilon'}{n+1}\, \mu^*\omega-\sum_{j\in J}\delta_j\,\zeta_j>0\,.
\end{eqnarray}
We notice that 
$$
0\;<\;\tau\;\le \;c^{\delta}_{j_0}\;=\;\frac{n-\delta_{j_0}}{h_{j_0}}\;<\;\frac{n}{n+1}\;<\;1\,,
$$
and that $\widehat\omega_{\delta}>0$ by the condition \eqref{Tie-bk-3}. The main up-shot of this choice of $\delta$ (usually called tie break) is that the condition \eqref{Tie-bk-2} implies $a^{\delta}_{j_1}=1$ and $a^{\delta}_j<1$ for all $j\in J:=\bar J\smallsetminus \{j_1\}$. This will allow us to restrict certain currents in an adequate way over the hypersurface $Z:=Z_{j_1}$. 
\\
\\
{\bf (C) Restriction and weak limit extraction.} 
\\
Let $\kappa:=\{\widehat\omega_{\delta}\}$ be the required K\"{a}hler class, let $I:=\{j\in J\,:\,a_j>0\}$ and let define the class $\alpha=\{[\Delta]\}+\kappa$  with respect to this data as in the statement of the lemma \ref{Sokur-const}.
We consider now a family of closed and real $(1,1)$-currents with analytic singularities 
$$
(R^{\varepsilon}_{\Theta})_{\varepsilon\in (0,1)}\in \{R_{\Theta}\}\,,
\qquad
R^{\varepsilon}_{\Theta}\;\ge\; -\varepsilon \,\widehat \omega\,,\qquad
\lambda(R^{\varepsilon}_{\Theta}\,,Z_j)=0\,,\quad \forall j\in \bar J\,.
$$
%with
%$$
%\lambda(T_{\varepsilon},y)\uparrow\lambda(\mu^{*}\Theta_{K_X+L},y)\,,
%\quad\mbox{as}\quad\varepsilon\downarrow 0\,,
%$$
obtained by regularising the current $R_{\Theta}$. 
We infer in particular the existence of the restriction
\begin{eqnarray}\label{lost-pos}
R^{\varepsilon}_{\Theta\,{|Z}}\ge -\varepsilon\,\widehat\omega_{{|Z}}\,.
\end{eqnarray}
The fact that 
$$
(R^{\varepsilon}_{\Theta\,{|Z}})_{\varepsilon\in (0,1)}\in \{R_{\Theta}\}_{{|Z}}\,,
$$ 
combined with the weak compactness of the mass imply the existence of a sequence $(\varepsilon_l)_l\subset(0,1) $ , $\varepsilon_l\downarrow 0$ as $l\rightarrow +\infty$ and a closed positive current 
$$
\Xi\in \{R_{\Theta}\}_{{|Z}}\,,
$$
such that $R^{\varepsilon_l}_{\Theta\,{|Z}}+\varepsilon_l\,\widehat\omega_{{|Z}}\longrightarrow \Xi$ weakly as $l\rightarrow +\infty$. Then the semi-continuity of the Lelong numbers implies
\begin{eqnarray}\label{sc-Lel-nr}
\lambda(\Xi, V_j)\;\ge\;\lambda(R^{\varepsilon_l}_{\Theta\,{|Z}}, V_j)\,,
\end{eqnarray}
for all $j\in J$ and all $l$. On the other hand the identity \eqref{key-idI} implies that the current
$$
\widehat\Theta_l\;:=\;\sum_{j\in J\smallsetminus I}-a_j[Z_j]\;+\;(\tau p+1)\,R^{\varepsilon_l}_{\Theta}\,,
$$
satisfies all the requirements in the statement of the lemma \ref{Sokur-const}. Moreover by restricting the identity \eqref{key-idI} to $Z$ we infer
$$
\widehat\Theta_Z\;:=\;\sum_{j\in J\smallsetminus I}-a_j[V_j]\;+\;(\tau p+1)\,\Xi\;\in \; c_1(K_Z)\;+\;\alpha_{{|Z}}\,.
$$
We decompose the current $\Xi$ as
$$
\Xi\;=\;\sum_{j\in I}\lambda(\Xi, V_j)\,[V_j]\;+\;R_{\Xi}\,,
$$
and we rewrite the closed positive current $\widehat\Theta_Z$ as
\begin{eqnarray}\label{key-rest-id-I}
\widehat\Theta_Z\;=\;\sum_{j\in I}\xi_j\,[V_j]\;+\;\widehat R_{\Xi}\;\in \; c_1(K_Z)\;+\;\alpha_{{|Z}}\,,
\end{eqnarray}
with $\xi_j\;:=\;(\tau p+1)\lambda(\Xi, V_j)$ and with
$$
\widehat R_{\Xi}\;:=\;\sum_{j\in J\smallsetminus I}-a_j[V_j]\;+\;(\tau p+1)\,R_{\Xi}\,.
$$
{\bf (D) Application of the non-vanishing in dimension one less.} 
\\
Let $\mu_j:=\xi_j-a_j$ for all $j\in I$. The identity \eqref{key-rest-id-I} implies
\begin{eqnarray}\label{key-rest-id-II}
\sum_{j\in I}\mu_j\, [V_j]\;+\;\widehat R_{\Xi}\;\in \; c_1(K_Z)\;+\;\kappa_{{|Z}}\,.
\end{eqnarray}
Let $I_+:=\{j\in I\,:\,\mu_j>0\}$, let $I_-:=I\smallsetminus I_+$ and let define the $(1,1)$-cohomology class
$$
\alpha_Z\;:=\;\sum_{j\in I_-}-\mu_j\, \{[V_j]\}\;+\;\kappa_{{|Z}}\,.
$$
We observe that $0\le -\mu_j\le a_j<1$ for all $j\in I_-$ and that the $(1,1)$-cohomology class $c_1(K_Z)+\alpha_Z$ is pseudoeffective
since
$$
0\;\le\;\sum_{j\in I_+}\mu_j\, [V_j]\;+\;\widehat R_{\Xi}\;\in \; c_1(K_Z)\;+\;\alpha_Z\,,
$$
by the identity \eqref{key-rest-id-II}. Thus the non-vanishing assumption in dimension on less implies the existence of an effective $\R$-divisor $\Gamma$ over $Z$ such that $[\Gamma]\in c_1(K_Z)+\alpha_Z$. On the other hand the decomposition
$$
\widehat\alpha_{{|Z}}\;=\;\alpha_Z\;+\;\sum_{j\in I_-}\xi_j\,\{[V_j]\}\;+\;\sum_{j\in I_+} a_j\,\{[V_j]\}\,,
$$
implies that the effective $\R$-divisor 
$$
G\;:=\;\Gamma\;+\;\sum_{j\in I_-}\xi_j\,V_j\;+\;\sum_{j\in I_+} a_j\,V_j\,,
$$
satisfies $[G]\in c_1(K_Z)+\widehat\alpha_{{|Z}}$. It also satisfies the vanishing condition \eqref{NV-Vansh-cnd} since
$$
\xi_j\;\ge\;(\tau p+1)\lambda(R^{\varepsilon_l}_{\Theta\,{|Z}}, V_j)\;=\;\lambda(\widehat\Theta_{l |Z},V_j)\,,
$$
for all $j\in I_-$ and $l$, by the inequality \ref{sc-Lel-nr}.
\\
\\
{\bf (E) End of the proof of the Shokurov's lemma.} 
\\
The identity \eqref{key-idII} decomposes as 
\eqref{Zot-Shok} with respect to $r:=\tau p+1$ and
$$
F\;:=\;\sum_{j\in I\cup \{j_1\}}(r\lambda_j+e_j)\,Z_j\;+\;\sum_{j\in J\smallsetminus I}[\tau h_j+(1-\tau)\theta_j+\delta_j]\,Z_j\,.
$$
\hfill$\Box$
%%%%%%%%%%%%%%%%%%%%%%%%%%%%%%%%%%%%%%%%%%%%%%%%%%%%%%%%%%%%%%%%%%%%%%%%%%%%%%%%%%%%%%%%%%%%%%%%%%%%%%%%%%%%%%%%%%%%%%%%%%%%%%%%%%
%%%%%%%%%%%%%%%%%%%%%%%%%%%%%%%%%%%%%%%%%%%%%%%%%%%%%%%%%%%%%%%%%%%%%%%%%%%%%%%%%%%%%%%%%%%%%%%%%%%%%%%%%%%%%%%%%%%%%%%%%%%%%%%%%%
%%%%%%%%%%%%%%%%%%%%%%%%%%%%%%%%%%%%%%%%%%%%%%%%%%%%%%%%%%%%%%%%%%%%%%%%%%%%%%%%%%%%%%%%%%%%%%%%%%%%%%%%%%%%%%%%%%%%%%%%%%%%%%%%%%
%%%%%%%%%%%%%%%%%%%%%%%%%%%%%%%%%%%%%%%%%%%%%%%%%%%%%%%%%%%%%%%%%%%%%%%%%%%%%%%%%%%%%%%%%%%%%%%%%%%%%%%%%%%%%%%%%%%%%%%%%%%%%%%%%%
%%%%%%%%%%%%%%%%%%%%%%%%%%%%%%%%%%%%%%%%%%%%%%%%%%%%%%%%%%%%%%%%%%%%%%%%%%%%%%%%%%%%%%%%%%%%%%%%%%%%%%%%%%%%%%%%%%%%%%%%%%%%%%%%%%
%%%%%%%%%%%%%%%%%%%%%%%%%%%%%%%%%%%%%%%%%%%%%%%%%%%%%%%%%%%%%%%%%%%%%%%%%%%%%%%%%%%%%%%%%%%%%%%%%%%%%%%%%%%%%%%%%%%%%%%%%%%%%%%%%%
%%%%%%%%%%%%%%%%%%%%%%%%%%%%%%%%%%%%%%%%%%%%%%%%%%%%%%%%%%%%%%%%%%%%%%%%%%%%%%%%%%%%%%%%%%%%%%%%%%%%%%%%%%%%%%%%%%%%%%%%%%%%%%%%%%
\subsection{Diophantine approximation}
The following lemma follows from quite elementary facts from linear algebra and diophantine approximation. The elementary details are left to the reader.
\begin{lem}\label{Diof-aprox}
Let
\\
$\bullet$  $\Delta=\sum_{j\in I}a_jZ_j$ be an effective $\Q$-divisor over a complex projective manifold $\widehat X$ such that $\lfloor\Delta\rfloor=0$,
\\
$\bullet$ $Z\subset \widehat X$ be a smooth and irreducible hypersurface which is not a component of the divisor $\Delta$ and $V_j:=Z\cap Z_j$ are distinct and irreducible for all $j\in I$,
\\
$\bullet$  $\widehat\alpha:=\{[\Delta]\}+\kappa$ be a $(1,1)$-cohomology class, with $\kappa \in \NS_{{\R}}(\widehat X)$ a K\"{a}hler class such that there exist an effective $\R$-divisor $G$ over $Z$ with the properties
$$
[G]\in c_1(K_Z)+\widehat\alpha_{{|Z}}\,,
$$
and
$$
G\;- \;\sum_{t\in I_-}\xi_t\,V_t\;-\;\sum_{t\in I_+} a_t\,V_t\;\ge \;0\,,\qquad 
\xi_t\in \R_{\ge 0}\,,\qquad
I=I_{+}\amalg I_{-}\,.
$$ 
%with $\xi_t\in \R_{\ge 0}$ and $I=I_{+}\amalg I_{-}$.
For any norm  $\|\cdot \|$ on the finite dimensional vector space $\NS_{{\R}}(\widehat X)$ there exist a constant $C>0$, a sequence $(m_p)_{p\in \N_{>0}}\subset \N_{>0}$ and finite families $(A^{j,p})_{j=1}^{N_p}$ of ample $\Q$-line bundles over $\widehat X$
such that for all $p\in \N_{>0}$,
\\
$\blacktriangleright$ $m_p\Delta$ is integral, $m_pA^{j,p}$ is a line bundle and
$$
\left\|c_1(A^{j,p})-\kappa\right\| \;\le \;\frac{C}{\,p\,m_p}\,,
$$
for all $j=1,...,N_p$.
\\
$\blacktriangleright$ there exist $r_{j,p}\in \R_{>0}$ such that 
$$
\kappa\;=\;\sum_{j=1}^{N_p}r_{j,p}\,c_1(A^{j,p})\,,\qquad\mbox{and}\qquad \sum_{j=1}^{N_p}r_{j,p}=1\,,
$$
\\
$\blacktriangleright$ if we set ${\cal L}^{j,p}:={\cal O}_{\widehat X}(\Delta)+A^{j,p}$ for all $j=1,...,N_p$, there exist an effective $\R$-divisor 
$$
[G^{j,p}]\in c_1(K_Z+{\cal L}^{j,p})\,,
$$
with the vanishing property
$$
G^{j,p}\;- \;\sum_{t\in I_-}\xi_t\,V_t\;-\;\sum_{t\in I_+} a_t\,V_t\;\ge \;-\;\frac{C}{m_p p}\,\sum_{t\in I_-}V_t\,.
$$
\end{lem}
%%%%%%%%%%%%%%%%%%%%%%%%%%%%%%%%%%%%%%%%%%%%%%%%%%%%%%%%%%%%%%%%%%%%%%%%%%%%%%%%%%%%%%%%%%%%%%%%%%%%%%%%%%%%%%%%%%%%%%%%%%%%%%%%%%
%%%%%%%%%%%%%%%%%%%%%%%%%%%%%%%%%%%%%%%%%%%%%%%%%%%%%%%%%%%%%%%%%%%%%%%%%%%%%%%%%%%%%%%%%%%%%%%%%%%%%%%%%%%%%%%%%%%%%%%%%%%%%%%%%%
%%%%%%%%%%%%%%%%%%%%%%%%%%%%%%%%%%%%%%%%%%%%%%%%%%%%%%%%%%%%%%%%%%%%%%%%%%%%%%%%%%%%%%%%%%%%%%%%%%%%%%%%%%%%%%%%%%%%%%%%%%%%%%%%%%
%%%%%%%%%%%%%%%%%%%%%%%%%%%%%%%%%%%%%%%%%%%%%%%%%%%%%%%%%%%%%%%%%%%%%%%%%%%%%%%%%%%%%%%%%%%%%%%%%%%%%%%%%%%%%%%%%%%%%%%%%%%%%%%%%%
%%%%%%%%%%%%%%%%%%%%%%%%%%%%%%%%%%%%%%%%%%%%%%%%%%%%%%%%%%%%%%%%%%%%%%%%%%%%%%%%%%%%%%%%%%%%%%%%%%%%%%%%%%%%%%%%%%%%%%%%%%%%%%%%%%
%%%%%%%%%%%%%%%%%%%%%%%%%%%%%%%%%%%%%%%%%%%%%%%%%%%%%%%%%%%%%%%%%%%%%%%%%%%%%%%%%%%%%%%%%%%%%%%%%%%%%%%%%%%%%%%%%%%%%%%%%%%%%%%%%%
%%%%%%%%%%%%%%%%%%%%%%%%%%%%%%%%%%%%%%%%%%%%%%%%%%%%%%%%%%%%%%%%%%%%%%%%%%%%%%%%%%%%%%%%%%%%%%%%%%%%%%%%%%%%%%%%%%%%%%%%%%%%%%%%%%
\subsection{Shokurov's construction of sections}
We remind now a well known fact due to Shokurov (see for example \cite{Pa2}).
\begin{Claim}\label{Shok-section}
Let $L$ be a holomorphic $\Q$-line bundle over a polarised connected complex projective manifold $(Z,\omega)$ which admits a closed positive $(1,1)$-current $\theta\in c_1(L)$ such that $\theta\ge \varepsilon\omega$ for some $\varepsilon\in \R_{>0}$ and such that ${\cal I}(\theta)={\cal O}_Z$. If there exists an effective $\Q$-divisor $G$ over $Z$ such that $[G]\in c_1(K_Z+L)$ then 
$$
h^0(Z, m(K_Z+L))>0\,,
$$
for all $m\in \N_{>0}$ such that $mL$ is a holomorphic line bundle and $mG$ is integral. Moreover if there exists an effective and simple normal crossing $\R$-divisor $V$ over $Z$ such that $G-V\ge 0$, then there exists a non zero section 
\begin{eqnarray}\label{Shok-prt-Sec-div}
u\in H^0(Z, m(K_Z+L))\,,\qquad \div u\;-\;\lfloor (m-1)V\rfloor\;\ge \; 0\,.
\end{eqnarray}
\end{Claim}
$Proof$. There exist a flat hermitian line bundle $F$ and a non zero section 
$$
\sigma\in H^0(Z,m(K_Z+L)+F)\,,
$$ 
such that $mG=\div \sigma$. Set 
$$
{\cal L}:=(m-1)(K_Z+L)+L\,,
$$ 
and observe the obvious identity $m(K_Z+L)=K_Z+{\cal L}$. We define the current
$$
\theta_{G}:=(m-1)[G]+\theta\in c_1({\cal L})=c_1({\cal L}+F)\,,\quad \theta_{G}\ge \varepsilon\omega\,,
$$
and we observe that $\sigma\in H^0(Z\,,\,{\cal S}_F)$, with
\begin{eqnarray*}%\label{Rat-sect-res}
{\cal S}_F\;:=\;{\cal S}\otimes_{{\cal O}_Z}{\cal O}_Z(F) \,,\qquad  {\cal S}\;:=\;  {\cal O}_Z(K_Z+{\cal L})\otimes_{{\cal O}_Z}{\cal I}(\theta_G)\,.
\end{eqnarray*}
In fact let $h$ be a smooth hermitian metric over 
$$
m(K_Z+L)+F\;=\;K_Z+{\cal L}+F\,,
$$ 
and let $\gamma_h\in c_1(K_Z+{\cal L})$ be its normalised curvature form. Let $\alpha\in c_1(L)$ smooth and let write $\theta=\alpha+\frac{i}{2\pi}\partial\bar\partial \varphi_{\theta}$. The Lelong-Poincar\'{e} formula implies
$$
\theta_G\;=\;\frac{m-1}{m}\,\gamma_h+\alpha+\frac{i}{2\pi}\,\partial\bar\partial \varphi_G\,,\qquad 
\varphi_G:=\frac{m-1}{m}\log|\sigma|^2_h+\varphi_{\theta}\,.
$$
Then 
$$
\int_Z|\sigma|^2_h\,e^{-\varphi_G}=\int_Z|\sigma|^{2/m}_h\,e^{-\varphi_{\theta}}\le C\int_Ze^{-\varphi_{\theta}}<+\infty\,,
$$
implies $\sigma\in H^0(Z\,,\,{\cal S}_F)$. 
By applying the  Kawamata-Viehweg-Nadel vanishing theorem to the line bundles ${\cal L}$ and ${\cal L}+F$ we infer
$$
h^q(Z\,,\,{\cal S})\;=\;h^q(Z\,,\,{\cal S}_F)\;=\;0\,,\qquad \forall q>0\,.
$$
Thus 
$$
h^0(Z\,,\,{\cal S})=\chi(Z\,,\,{\cal S})=\chi(Z\,,\,{\cal S}_F)=h^0(Z\,,\,{\cal S}_F)>0\,,
$$
%On the other hand the fact that $\Ch(F)=1$ implies $\Ch({\cal S})=\Ch({\cal S}_F)$ by basic facts about Chern caracters. Then using the Hirzebruch-Riemann-Roch formula we infer
%\begin{eqnarray*}
%h^0(Z\,,\,{\cal S})&=&\chi(Z\,,\,{\cal S})=\int_Z\Td(Z)\cdot \Ch({\cal S})
%\\
%\\
%&=&
%\chi(Z\,,\,{\cal S}_F)=h^0(Z\,,\,{\cal S}_F)>0\,. 
%\end{eqnarray*}
%\begin{eqnarray*}h^0(Z\,,\,{\cal S})&=&\chi(Z\,,\,{\cal S})\;=\;\chi(Z\,,\,{\cal S}\otimes_{{\cal C}^0_Z}{\cal C}^0_Z )\\\\&=&
%\chi(Z\,,\,{\cal S}_F\otimes_{{\cal C}^0_Z}{\cal C}^0_Z)=\chi(Z\,,\,{\cal S}_F)=h^0(Z\,,\,{\cal S}_F)>0\,,\end{eqnarray*}
since $F$ is topologically trivial.
Moreover the inclusion 
$$
{\cal I}(\theta_G)\subset {\cal I}((m-1)[V])={\cal O}_Z(-\lfloor (m-1)V\rfloor)\,,
$$
implies the existence of the required section $u$.
\hfill$\Box$
\\
\\
{\bf Conclusion}.  It seem clear at this point that the magnitude of the vanishing error of type  \eqref{Shok-prt-Sec-div} produced by a combination of diophantine approximation with   Shokurov's construction of sections is much bigger than the magnitude of the vanishing error allowed by the extension condition  \eqref{prt-Sec-div} of corollary \ref{Pert-Hk-Mk-lem}. To be more precise we want to apply claim \ref{Shok-section} in the setting of lemma \ref{Sokur-const} after application of the diophantine approximation of lemma \ref{Diof-aprox}. Therefore we want to apply the claim \ref{Shok-section} with $L:= {\cal L}^{j,p}$ and with
$$
V\;=\;\sum_{t\in I}\beta_t\,V_t\;:=\;\sum_{t\in I_-}\left(\xi_t\;+\;\frac{C}{ m_p p}\right)V_t\;+\;\sum_{t\in I_+} a_t\,V_t\,.
$$
The trivial identity 
$$
\beta_t\;-\;\frac{1}{m_p}\,\lfloor (m_p-1)\beta_t\rfloor\;=\;\frac{1}{m_p}\,\beta_t\;+\;\frac{1}{ m_p}\Big[(m_p-1)\beta_t\;-\;\lfloor (m_p-1)\beta_t\rfloor\Big]\,,
$$
shows that the error 
$$
\frac{1}{m_p}\,\beta_t\;\le\;\beta_t\;-\;\frac{1}{m_p}\,\lfloor (m_p-1)\beta_t\rfloor\,,
$$
does not allow to apply the condition  \eqref{prt-Sec-div} in corollary \ref{Pert-Hk-Mk-lem} to the section $u$ produced by claim \ref{Shok-section}. We observe furthermore that if one let $m\rightarrow +\infty$ the condition \eqref{err-pos-curr} in corollary \ref{Pert-Hk-Mk-lem} force a priori a blow-up of the generic Lelong numbers $\xi_t=\lambda(\Theta_{|Z}\,,V_t)$ and therefore a blow-up of $\beta_t$, $t\in I_-$. We observe finally that the situation $I_+=\emptyset$ my occur in lemma \ref{Sokur-const}.  
\\
\\
{\bf Acknowledgments.} The author is grateful to Professors Adrien Dubouloz and Gabriele La Nave for useful conversations.

\vspace{1cm}
\noindent
Nefton Pali
\\
Universit\'{e} Paris Sud, D\'epartement de Math\'ematiques 
\\
B\^{a}timent 425 F91405 Orsay, France
\\
E-mail: \textit{nefton.pali@math.u-psud.fr}
\end{document}